\documentclass[12pt]{amsart}
\usepackage{graphicx}
\usepackage{amsmath}
\usepackage{amsfonts}
\usepackage{amssymb}
\usepackage{amscd}
\usepackage{amsthm}

\newtheorem{theorem}{Theorem}

\newtheorem{lemma}[theorem]{Lemma}
\newtheorem{proposition}[theorem]{Proposition}

\newtheorem{remark}[theorem]{Remark}

\def\x{\xi}

\def\b{\beta}

\def\f{\psi}
\def\S{\mathbf{S}}
\def\Y{\mathcal{Y}}
\def\MCG{\textup{MCG}}
\def\PMCG{\textup{PMCG}}
\def\H{\textup{H}}
\def\F{\textup{T}}
\def\V{V}
\def\E{\mathcal{E}}

\date{}

\begin{document}

\title[Extending homeomorphisms from 2-punctured surfaces]{Extending homeomorphisms from 2-punctured surfaces to handlebodies}

\author{Alessia Cattabriga \and  Michele Mulazzani}

\begin{abstract}
Let $\H_g$ be a genus $g$ handlebody, and $\F_g=\partial\H_g$ a
closed connected orientable surface.  In this paper we find a
finite set of generators for $\E_2^g$, the subgroup of
$\PMCG_2(\F_g)$ consisting of the isotopy classes of
homeomorphisms of $\F_g$ which admit an extension to the
handlebody keeping  a properly embedded trivial arc fixed. This
subgroup turns out to be important for the study of knots in
closed 3-manifolds via $(g,1)$-decomposition. In fact,  the knots
represented by the isotopy classes belonging to the same left
cosets of $\E_2^g$ in $\PMCG_2(\F_g)$ are equivalent.

\bigskip

\noindent{{\it Mathematics Subject
Classification 2000:} Primary  20F38; Secondary  57M25, 57N10.}

\end{abstract}

\maketitle

\section{Introduction} \label{intro}

The study of knots and links in closed 3-manifolds via Heegaard
decompositions was started by Doll in \cite{Do}, where the notion
of $(g,b)$-decompositions of links is introduced. In particular,
the case of $(1,1)$-decompositions of knots in lens spaces and in
the 3-sphere has been extensively investigated in several papers
(see for example \cite{CM1,CM2,CK,GHS,GMM,Ha,Sa}). More generally,
any knot $K$ in a 3-manifold $M$ of Heegaard genus $h$ admits a
$(g,1)$-decomposition for a suitable $g\geq h$. This means that
there exists a genus $g$ Heegaard splitting $(\H_g,\H'_g)$ of $M$
such that $\F_g=\partial\H_g=\partial\H'_g$ is transverse to $K$
and splits it in two arcs, $A\subset\H_g$ and $A'\subset\H'_g$,
parallel to $\F_g$ (for details see \cite{CMV}). As a consequence,
the knot $K$ can be represented, up to equivalence, by an element
of $\PMCG_2(\F_g)$, the pure mapping class group of the twice
punctured surface of genus $g$. Moreover, if a homeomorphism
$\psi:(\partial\H_g,\partial A)\to(\partial\H_g,\partial A)$
admits an extension to the handlebody $\H_g$ fixing $A$, the knot
$K_{\psi}$, represented by the isotopy class of $\psi$, is the
trivial knot in the connected sum of $g$-copies of
$\S^1\times\S^2$. Therefore, if we denote with $\E_2^g$ the
subgroup of $\PMCG_2(\F_g)$, consisting of the elements which
admit an extension to $\H_g$ fixing $A$, and $\mathcal{K}_{g,1}$
is the set of all the knots admitting a $(g,1)$-decomposition,
there exists a surjective map
$$\Theta_g:\PMCG_2(\F_g)/\E_2^g\to\mathcal{K}_{g,1},$$ where
$\PMCG_2(\F_g)/\E_2^g$ denotes the set of left cosets of  $\E_2^g$
in $\PMCG_2(\F_g)$.

The aim of this paper is to find a finite set of  generators for
$\E_2^g$. In Section \ref{generators}, we introduce notations,
describe some significant elements of  $\E_2^g$, which have an
explicit topological meaning, and recall some known facts. The
main result of the paper is given in  Section \ref{main result},
where a finite set of generators for $\E_2^g$ is obtained, for
each $g\geq 1$.

\section{Preliminaries and notations}
\label{generators}

Given a closed orientable connected surface $\F_g$ of genus $g\geq
0$ and $n$ distinguished points  $P_1,\ldots,P_n$ on it, called
\textit{punctures}, we denote with $\mathcal{H}_n(\F_g)$ (resp.
$\mathcal{F}_n(\F_g)$) the group of orientation-preserving
homeomorphisms $h:\F_g\to \F_g$ such that
$h(\{P_1,\ldots,P_n\})=\{P_1,\ldots,P_n\}$ (resp. $h(P_i)=P_i$ for
$i=1,\ldots,n$). The \textit{ $n$- punctured mapping class group}
$\MCG_n(\F_g)$ (resp. the \textit{$n$-punctured pure mapping class
group} $\PMCG_n(\F_g)$) of $\F_g$ is  the group of the isotopy
classes of elements of $\mathcal{H}_n(\F_g)$ (resp.
$\mathcal{F}_n(\F_g)$).  Clearly, neither $\MCG_n(\F_g)$ nor
$\PMCG_n(\F_g)$ depend on the choice of the punctures on $\F_g$.
There is an exact sequence
$$1\to \PMCG_n(\F_g) \to \MCG_n(\F_g)\to \Sigma_n \to 1,$$
where  $\Sigma_n$ is the symmetric group of degree $n$. Moreover,
we denote with $\MCG(\F_g)$ the mapping class group of $\F_g$ with
no punctures.

Let $\H_g$ be  an orientable handlebody of genus $g$ and denote
with $A\subset \H_g$ a properly embedded trivial arc (i.e., there
exists a disk $D\subset \H_g$ with $A\cap D=A\cap\partial D=A$ and
$\partial D-A\subset\partial \H_g$). The boundary  of $\H_g$ is a
closed orientable connected surface of genus $g$ with two
distinguished punctures $P_1$ and $P_2$ on it, which are the
endpoints of $A$. In this article we are interested in determining
generators for  the subgroup $\E_2^g$ of $\PMCG_2(\F_g)$,
consisting of the elements which admit an extension to $\H_g$
fixing $A$.

\begin{figure}[ht]
\begin{center}
\includegraphics*[totalheight=8cm]{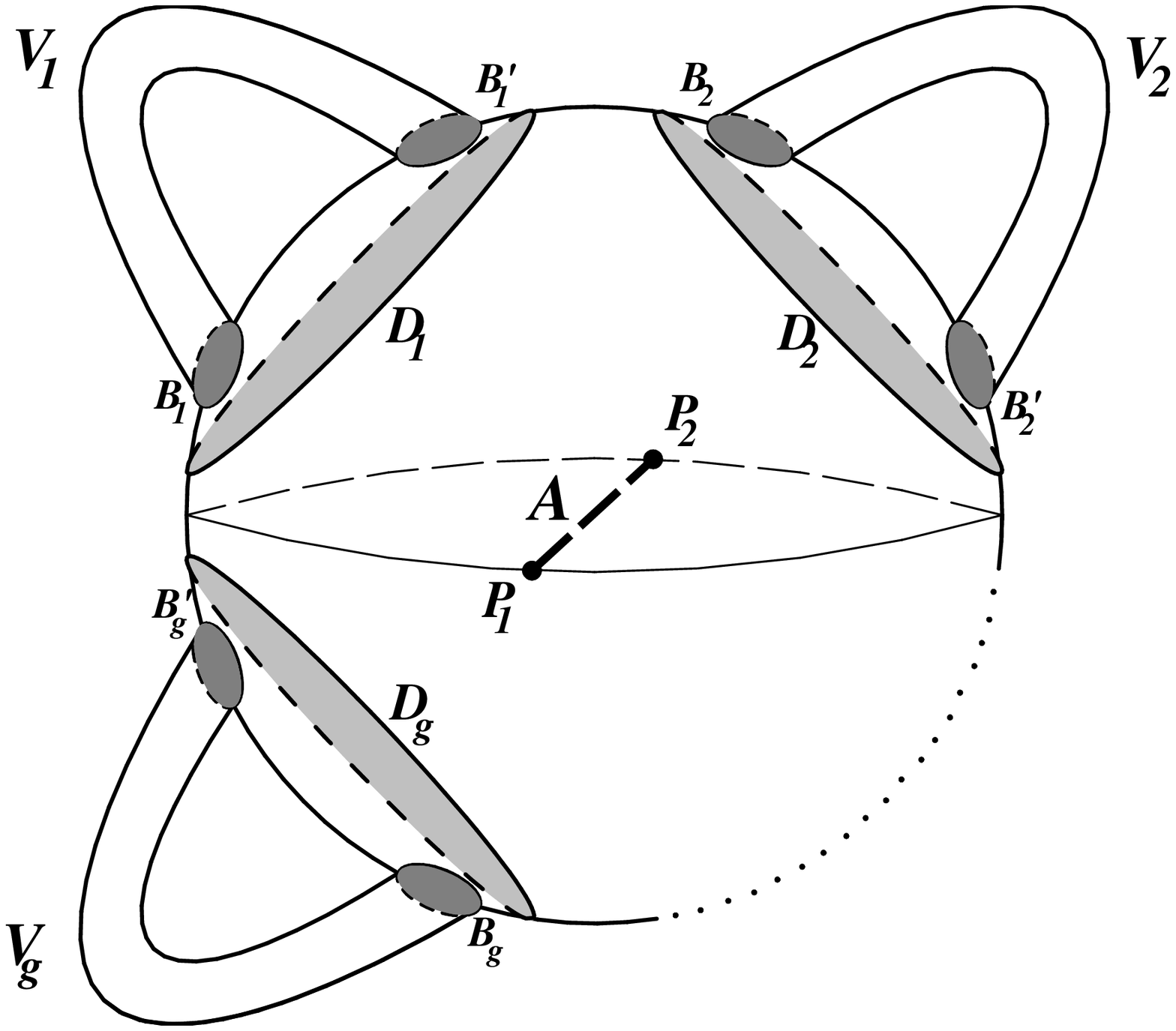}
\end{center}
\caption{A model for $\H_g$.} \label{Fig1}
\end{figure}

The model for a genus $g$ handlebody $\H_g$ will be the one
depicted in Figure \ref{Fig1}, and $A\subset \H_g$ will be the
properly embedded trivial arc dashed  in the figure.

In order to simplify notation, we will use the same symbol to
denote  a homeomorphism and its isotopy class. Moreover,
$\overline{\psi}$ will denote an extension to $\H_g$ of a
homeomorphism $\psi$ of $\F_g$. The right-handed Dehn twist along
a simple closed curve $e$ will be denoted with $t_{e}$, while
$s_{P,e}$ will denote the spin of $P$ about $e$ (see \cite{Bi}).

Now we briefly recall some homeomorphisms of $\F_g$ that admit an
extension to $\H_g$ keeping $A$ fixed (for details see \cite{Su}).

\begin{description}
\item[Rotation] let  $\overline{\rho}$  be the clockwise
$\frac{2\pi}{g}$-rotation  of $\H_g$  around the  axis containing
the punctures $P_1$ and $P_2$. The restriction   of this
homeomorphism to $\F_g$ will be denoted with $\rho$.

\item[Semitwist of a knob] let $D_i$ be the   2-cell depicted in
Figure~\ref{Fig1} and let $\delta_i=\partial D_i$, for
$i=1,\ldots,g$. Then $D_i$ cuts away from $\H_g$ a solid torus
$K_i$, containing the $i$-th handle $\V_i$, that will be called
the $i$-th knob. We denote with $\overline{\omega}_i$ the
homeomorphism of $\H_g$ obtained by cutting $\H_g$ along $D_i$,
giving  a semitwist  counterclockwise to $K_i$ and gluing it back,
for $i=1,\ldots, g$. Moreover, we set
$\omega_i=(\overline{\omega}_i)_{|_{\F_g}}$. Notice that
$\omega_i^2$ is isotopic to the Dehn twist along $\delta_i$.

\item[Twisting a handle]  Referring to Figure \ref{Fig1}, let
$\b_i=\partial B_i$, for $i=1,\ldots, g$. We denote with $\tau_i$
the right-handed Dehn twist along $\b_i$. Note that $\tau_i$
admits an extension to $\H_g$ whose effect is to give a complete
twist to the $i$-th handle, for $i=1,\ldots, g$.

\item [Slides] Referring to Figure \ref{Fig1}, let $Z_{i}$ (resp.
$Z'_i$) be the center of the properly embedded meridian disc
$B_i\subset\H_g$ (resp. $B'_i\subset\H_g$). Moreover, denote with
$\H_g^i$ the genus $(g-1)$ handlebody obtained  by removing the
$i$-th handle $\V_i$ from $\H_g$, and let $e$ (resp. $e'$) be a
simple closed curve on $\partial\H_g^i-\{P_1,P_2\}$ containing
$Z_i$ (resp. $Z'_i$), with $Z'_i \notin e$ (resp. $Z_i \notin
e'$). Such a curve will be called an $i$\textit{-loop} (resp. an
$i'$-\textit{loop}). Consider the spin $s_{Z_i,e}$ on
$\partial\H_g^i$. We can suppose, up to isotopy, that
$s_{Z_i,e}(B_i)=B_i$ and $s_{Z_i,e}(B'_i)=B'_i$.  So, extending
$s_{Z_i,e}$ by the identity on $\partial \V_i$, we obtain a
homeomorphism of $\F_g$ which will be denoted with $\sigma_{i,e}$
and called the \textit{slide} of $Z_i$ along $e$. In a completely
analogous way, given a $i'$-loop $e'$, we can define a slide of
$Z'_i$ along $e'$ and denote it with $\sigma'_{i',e'}$. Let
$\mathcal S=\{\sigma_{i,e}\,|\,e \textup{ is an } i \textup{-loop
}, i=1,\ldots,g\}\cup \{\sigma'_{i',e'}\,|\,e' \textup{ is an } i'
\textup{-loop }, i'=1,\ldots,g\}$.  Moreover, we set
$\theta_{ij}=\sigma_{i,e_{ij}}$, $\xi_{ij}=\sigma_{i,g_{ij}}$,
$\eta_{ik}=\sigma_{i,f_{ik}}$ (resp.
$\theta'_{ij}=\sigma'_{i',e'_{ij}}$,
$\xi'_{ij}=\sigma'_{i',g'_{ij}}$,
$\eta'_{ik}=\sigma'_{i',f'_{ik}}$), where
$e_{ij},g_{ij},f_{ik},e'_{ij},g'_{ij},f'_{ik}$ are the curves
depicted in Figure \ref{Fig2}, with $i,j=1,\ldots, g$, $i\neq j$,
and $k=1,2$.

\begin{figure}[ht]
\begin{center}
\includegraphics*[totalheight=7cm]{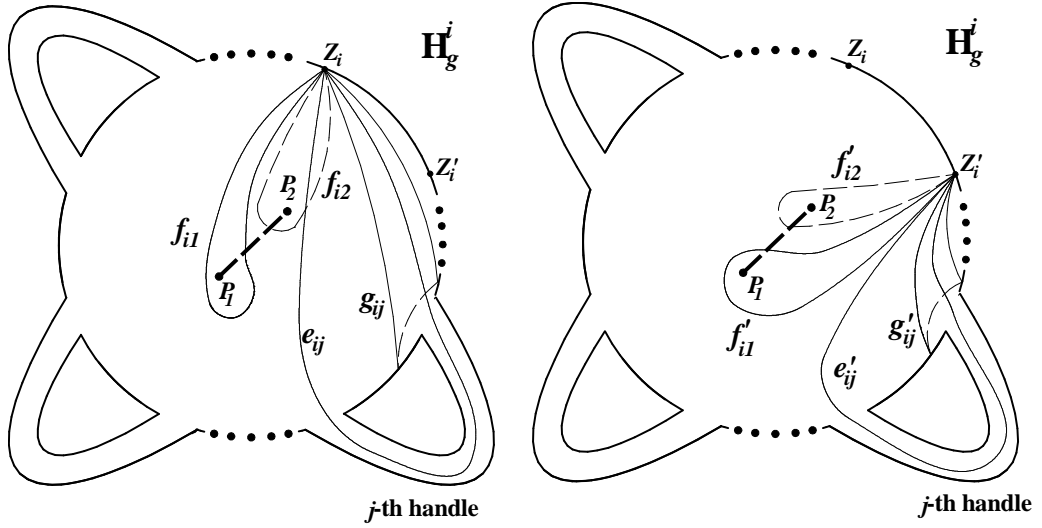}
\end{center}
\caption{The loops $e_{ij},g_{ij},f_{ik},e'_{ij},g'_{ij},f'_{ik}$,
on $\H^i_g$.} \label{Fig2}
\end{figure}

\item[Exchanging  two knobs] referring to Figure \ref{Fig4}, let
$C_i$ (resp. $C_j$) be the 2-cell on $\partial H_g^i$ (resp.
$\partial H_g^j$) bounded by $\nu_i$ (resp. $\nu_j$) not
containing $K_j$ (resp. $K_i$). Denote with $D$ a tubolar
neighborhood of $C_i\cup e\cup C_j$ on $\F_g$. Let $\rho'_{ij}$ be
the homeomorphism of $\F_g$ obtained  by exchanging $\partial K_i$
with $\partial K_j$ in such a way that $\rho'_{ij}$ twists
$C_i\cup e\cup C_j$ inside $D$ through $\pi$ radians in the
clockwise direction. So, $\rho'_{ij}(C_i)=C_j$,
$\rho'_{ij}(C_j)=C_i$ $\rho'_{ij}(e)=e$ and $\rho'_{ij}$ is the
identity outside a tubolar neighborhood of $K_i\cup e\cup K_j$. We
set $\rho_{ij}=\omega_j\omega_i\rho'_{ij}$, for $i,j=1,\ldots, g$
and $i\ne j$.
\end{description}
\begin{figure}[ht]
\begin{center}
\includegraphics*[totalheight=6cm]{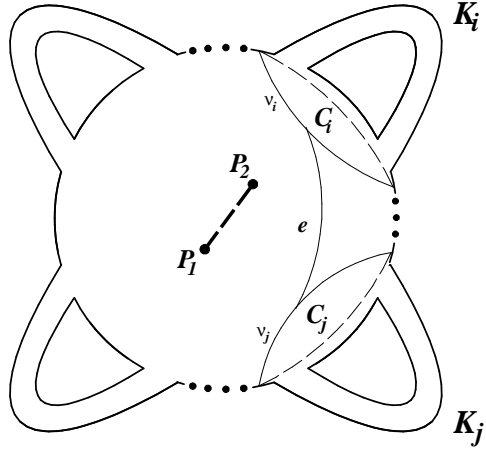}
\end{center}
\caption{Exchanging two knobs.} \label{Fig4}
\end{figure}

\begin{proposition}[\cite{Su}]\label{spin} Let $e_0,\ldots ,e_n$ (resp. $e'_0,\ldots ,e'_n$) be
$i$-loops (resp. $i'$-loops) such that $e_0$ (resp. $e'_0$) is
homotopic to $e_1\cdots e_n$ (resp. $e'_1\cdots e'_n$) on
$\partial \H_g^i-\{P_1,P_2,Z'_i\}$ rel $Z_i$ (resp. $\partial
\H_g^i-\{P_1,P_2,Z_i\}$ rel $Z'_i$). Then $\sigma_{i,e_0}$ (resp.
$\sigma'_{i',e'_0}$) is isotopic to
$\sigma_{i,e_1}\cdots\sigma_{i,e_n}$ (resp.
$\sigma'_{i',e'_1}\cdots\sigma'_{i',e'_n}$) modulo $\tau_i$.
\end{proposition}

\begin{proposition} [\cite{Su}] The subgroup $\mathcal E^g$ of $\MCG(\F_g)$, containing the elements of $\MCG(\F_g)$
that extend to $\H_g$  is generated by
$\rho,\rho_{12},\xi_{12},\theta_{12},\tau_1,\omega_1$.
\end{proposition}

\section{Main result}
\label{main result} In this paragraph we establish the main result
of the article.

\begin{figure}[ht]
\begin{center}
\includegraphics*[totalheight=6.5cm]{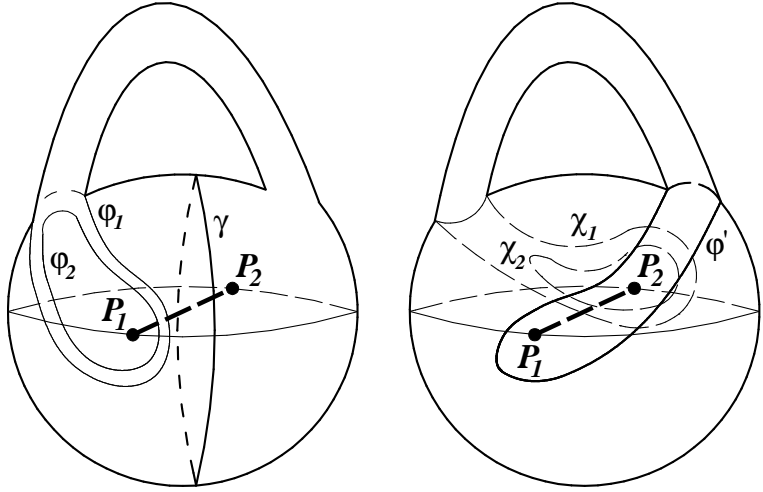}
\end{center}
\caption{The case of genus one.} \label{Fig5}
\end{figure}

\begin{theorem} \label{main}  The subgroup $\E^g_{2}$ of $\PMCG_2(\F_g)$
is generated by
$\rho,\rho_{12},\xi_{12},\theta_{12},\tau_1,\omega_1,\eta_{11},\eta_{12}$
if $g>1$, and $\tau_1,\omega_1,t_{\gamma}$ if $g=1$, where
$\gamma$ is the curve depicted in Figure \ref{Fig5}.
\end{theorem}

Before giving the proof we need some preparatory lemmas. Denote
with $\mathcal{Y}^g$ the subgroup of $\E_2^g$ generated by
$\rho_{ij},\xi'_{ij},\xi_{ij},\theta'_{ij},\theta_{ij},\tau_i,\omega_i,\eta_{i1},\eta_{i2},\eta'_{i1},\eta'_{i2}$,
with $i,j=1,\ldots ,g$ and $i\ne j$.

\begin{lemma} \label{slide} Let  $e$ (resp. $e'$) be  an $i$-loop (resp. $i'$-loop). Then
$\sigma_{i,e}$ (resp. $\sigma'_{i,e'}$)  is isotopic to a product
of $\x_{ij},\theta_{ij},\eta_{ik},\tau_{i}$ (resp.
$\x'_{ij},\theta'_{ij},\eta'_{ik},\tau_{i}$), for $j=1,\ldots,g$,
$k=1,2$ and $i\ne j$.
\end{lemma}

\begin{proof}
Since  $e_{ij},g_{ij},f_{ik}$, (resp. $e'_{ij},g'_{ij},f'_{ik}$)
with $j=1,\ldots,n$, $k=1,2$ and $i\ne j$, are a free set of
generators of $\pi_1(\partial \H_g^i-\{P_1,P_2,Z'_i\},Z_i)$,
(resp. $\pi_1(\partial \H_g^i-\{P_1,P_2,Z_i\},Z'_i)$) the lemma
follows from Proposition~\ref{spin}.
\end{proof}
\begin{lemma} \label{suzuki}
For each $f\in\E^g_2$ there exists an element $f_0\in \Y_g$ such
that $\overline{f}_0\overline{f}$ is the identity on the $i$-th
handle.
\end{lemma}
\begin{proof}
By Lemma \ref{slide}, $\mathcal{S}\subseteq\Y_g$. Therefore, using
the same arguments as in the proof of \cite[Lemma 4.4]{Su} we have
that $\overline{f}_0\overline{f}$ is the identity on a meridian
disk of the $i$-th handle and so, up to isotopy, on all the
handle.
\end{proof}

\begin{proposition}\label{braid}For $g\geq 1$ and $1\leq m<n$, consider the homomorphism
$j_{g,n,m}:\PMCG_n(\F_g)\rightarrow\PMCG_m(\F_g)$ induced by the
inclusion. Then $\ker j_{g,n,m}\cong \pi_1(F_{m,n-m}(\F_g))$,
where $F_{m,n-m}(\F_g)$ denotes the configuration space of $n-m$
points in $\F_g-\{P_1,\ldots P_m\}$.
\end{proposition}
\begin{proof} Let $P_1,\ldots,P_n$ be distinct points on
$\F_g$. As in \cite[Theorem 4.1]{Bi}, it is possible to prove that
the evaluation map $\epsilon: \mathcal{F}_m(\F_{g})\to
F_{m,n-m}(\F_g)$, given by
$\epsilon(\f)=(\f(P_{m+1}),\ldots,\f(P_{n}))$, is a locally
trivial fibering with fibre $\mathcal{F}_n(\F_{g})$, so the
following exact sequence holds:
$$\pi_1(F_{m,n-m}(\F_g),(P_{m+1},\ldots,P_n))\stackrel{d_{\#}}{\rightarrow}
\pi_0(\mathcal{F}_n(\F_{g}),id)
\stackrel{j_{\#}}{\rightarrow}\pi_0(\mathcal{F}_m(\F_{g}),id)\rightarrow
1.$$ Since $\pi_0(\mathcal{F}_k(\F_{g}),id)=\PMCG_k(\F_g)$ and
$j_{\#}$  is induced by the inclusion, it is enough to prove that
$\ker(d_{\#})=1$ to get the statement. Proceeding as in the proof
of \cite[Lemma 4.2.1]{Bi}, we obtain that $\ker(d_{\#})$ is
contained in the center of
$\pi_1(F_{m,n-m}(\F_g),(P_{m+1},\ldots,P_n))$. Now we will prove,
by induction on $k$,  that the center of $\pi_1(F_{m,k}(\F_g),*)$
is trivial. For $k=1$ we have
$\pi_1(F_{m,1}(\F_g),*)\cong\pi_1(\F_g-\{P_1,\ldots,P_m\},*)$,
which is a free group on $2g+m-1>1$ generators. Now let $k>1$. By
\cite[Theorem 1.2]{Bi} for each $1\leq r<k$ the projection
$\pi:F_{m,k}(\F_g)\to F_{m,r}(\F_g)$ is a locally trivial fibering
with fiber $F_{m+r,k-r}(\F_g)$. Since
$\pi_3(F_{m,1}(\F_g))=\pi_3(\F_g-\{P_1,\ldots,P_m\})=1$ and
$\pi_2(F_{m,1}(\F_g))=\pi_2(\F_g-\{P_1,\ldots,P_m\})=1$, from the
long exact sequence of the above fibration, for $r=1$, we get
$\pi_2(F_{m+1,k-1}(\F_g))\cong\pi_2(F_{m,k}(\F_g))$. So an
inductive argument shows that \hbox{$\pi_2(F_{m,k-1}(\F_g))\cong
\pi_2(F_{m+1,k-2}(\F_g))\cong\cdots\cong\pi_2(F_{m+k-2,1}(\F_g))=1$.}
Then, from the long exact sequence of the fibration in the case
$r=k-1$ we get
$$1\to\pi_1(F_{m+k-1,1}(\F_g))\stackrel{i_{\#}}{\rightarrow}\pi_1(F_{m,k}(\F_g))
\stackrel{\pi_{\#}}{\rightarrow}\pi_1(F_{m,k-1}(\F_g))\to 1.$$
Now, by the induction hypothesis, the center of
$\pi_1(F_{m,k-1}(\F_g))$ is trivial. So the center of
$\pi_1(F_{m,k}(\F_g))$ is contained in $\ker(\pi_{\#})\cong
i_{\#}( \pi_1(F_{m+k-1,1}(\F_g)))\cong \pi_1(F_{m+k-1,1}(\F_g))$,
which is centerless.
\end{proof}
\begin{remark} \label{ker}\textup{ As in \cite[pp. 158-160]{Bi}, a system of generators
of $\ker j_{g,n,m}$ can be constructed as follows. Let $E_n$ be a
set of generators of \hbox{$\pi_1(\F_g-\{P_1,\ldots P_m\},P_n)$},
such that for each $e\in E_n$, we have
$e\cap\{P_1,\ldots,P_n\}=\{P_n\}$. Moreover, for each
$i=m+1,\ldots,n-1$, denote with $E_i$ a set of generators of
$\pi_1(\F_g-\{P_1,\ldots P_m\},P_i)$ such that
$e\cap\{P_1,\ldots,P_n\}=\{P_i\}$ and each loop of $E_i$ is free
homotopic on $\F_g-\{P_1,\ldots,P_m\}$ to a loop in $E_n$. Then
$\ker j_{g,n,m}$  is generated by  $\{s_{i,e}\,|\,e\in E_i,
i=m+1,\ldots,n\}$.}
\end{remark}

Now we are ready to prove  Theorem \ref{main}.
\begin{proof}
Suppose that $\mathcal{Y}^g=\E_2^g$. From the relations in
\cite[Section 3]{Su}, and since
$\eta_{ik}=\rho^{i-1}\eta_{1k}\rho^{-(i-1)}$,
$\eta'_{i1}=\omega_i^{-1}\eta_{i1}\omega_i$ and
$\eta'_{i2}=\omega_i\eta_{i2}\omega_i^{-1}$, for $i=1,\ldots,g$,
the subgroup $\mathcal{Y}^g$ is generated by
$\rho,\rho_{12},\xi_{12},\theta_{12},\tau_1,\omega_1,\eta_{11},\eta_{12}$
if $g>1$, and by $\tau_1,\omega_1,\eta_{11},\eta_{12}$ if $g=1$.
This completes the proof when $g>1$. If  $g=1$,  referring to
Figure \ref{Fig5}, we have
$\eta_{11}=t_{\varphi_1}t_{\varphi_2}^{-1}= t_{\varphi_1}$ and
$\eta_{12}=t_{\chi_1}t_{\chi_2}^{-1}= t_{\chi_1}$. Moreover,
$\varphi_1\simeq\gamma$, $\chi_1\simeq\varphi'$ and
$\omega_1(\varphi_1)=\varphi'$.

Now we show, by induction on the genus, that
$\mathcal{Y}^g=\E_2^g$. Let $\widehat{\E}_2^g$ be the subgroup of
$\E_2^g$ containing the isotopy classes of the homeomorphisms that
are the identity on the boundary of the $g$-th handle $\V_g$. By
Lemma \ref{suzuki}, up to composing with an element of
$\mathcal{Y}^g$, we can suppose that each  $f\in\E_2^g$ belongs to
$\widehat{\E}_2^g$. So it is enough to show  that
$\widehat{\E}_2^g\subseteq\mathcal{Y}^g$.

Let $f:\F_g\to\F_g$ be a homeomorphism fixing $\partial \V_g$
pointwise and whose isotopy class belongs to $\widehat{\E}_2^g$.
By cutting out the $g$-th handle, and capping  with the two disks
$B_g$ and $B'_g$,  we can identify $f$ with a homeomorphism $f'$
of $\F_{g-1}$, such that $f'_{|_{B_g\cup B'_g}}=\textup{Id}$ and
$f'=f$ on $\F_g-(B_g\cup B'_g)$. Moreover, by shrinking $B_g$ and
$B'_g$ to their centers $Z_g$ and $Z'_g$, the map $f'$ becomes a
map $\widetilde{f}$ of $\F_{g-1}$ fixing $Z_g$ and $Z'_g$. In
order to simplify the notation, we set $P_3=Z_g$ and $P_4=Z'_g$.
Obviously, $\widetilde{f}$ extends to   $\H_{g-1}$ fixing $A$.
Note that isotopies of $\widetilde{f}$ rel $\{P_1,P_2,P_3,P_4\}$
correspond to isotopies of $f$ rel $\{P_1,P_2\}$ and the isotopy
class of $\widetilde{f}$ is trivial in $\PMCG_4(\F_{g-1})$ if and
only if the isotopy class of $f$ belongs to the subgroup generated
by $\tau_g$ in $\PMCG_2(\F_{g})$.

Let $\widetilde{\E}^{g-1}$ be the subgroup of $\PMCG_4(\F_{g-1})$
consisting of  the elements which extend to the handlebody keeping
$A$ fixed. Note that for each $g\in\widetilde{\E}^{g-1}$ we have
$g=\widetilde{f}$, for a suitable $f$.
\begin{figure}[ht]
\begin{center}
\includegraphics*[totalheight=3.5cm]{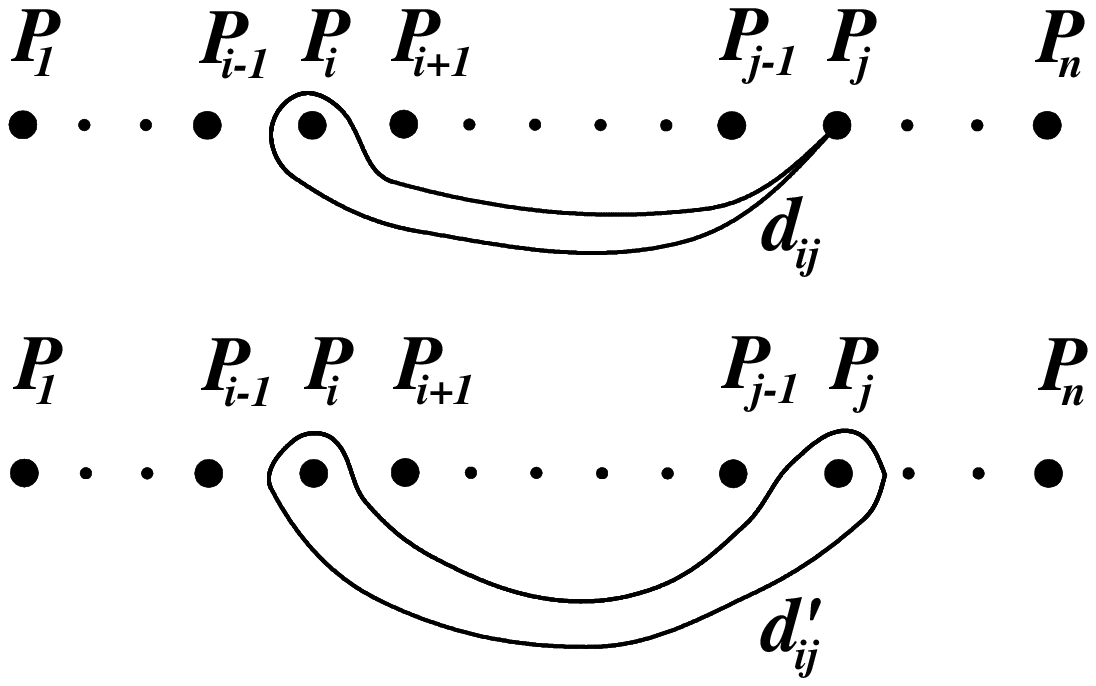}
\end{center}
\caption{Spin $s_{ij}$ on $\S^2$.} \label{Fig3}
\end{figure}
Let $g=1$.  We have that
$\PMCG_4(\S^2)=F_{4,0}/\textup{center}(F_{4,0})$ and a set of
generators for $\PMCG_4(\S^2)$ is  given by
$\{s_{ij}\,|\,i,j=1,\ldots,4;i< j\}$, where $s_{ij}$ is the spin
of $P_j$ along the curve $d_{ij}$ depicted in Figure \ref{Fig3}
(see \cite[Lemma 1.8.2, Theorem 4.2]{Bi}), which is equal to the
Dehn twist along the curve $d'_{ij}$. Obviously, $t_{d'_{ij}}$
extends to the 3-ball fixing $A$, and therefore
$\PMCG_4(\S^2)=\widetilde{\E}^0$. It is easy to see that
$s_{14}=\widetilde{\eta'_{11}}$, $s_{13}=\widetilde{\eta_{11}}$,
$s_{24}=\widetilde{\eta'_{12}}$, $s_{23}=\widetilde{\eta_{12}}$
and, since $d'_{12}$ is isotopic to $d'_{34}$ on $\S^2$ rel
$\{P_1,P_2,P_3,P_4\}$, we have $s_{12}=s_{34}=
\widetilde{t_{\delta_1}}=\widetilde{\omega_1^2}$. So the isotopy
class of $f$ is generated by
$\eta'_{11},\eta_{11},\eta'_{12},\eta_{12},\omega_1^2,\tau_1$. As
a consequence $\widehat{\E}_2^1\subseteq\mathcal{Y}^1$.

Now let $g>1$. The homomorphism
$j_{g-1,4,2}:\PMCG_4(\F_{g-1})\rightarrow\PMCG_2(\F_{g-1})$ of
Proposition \ref{braid} restricts to a homomorphism
\hbox{$j=(j_{g-1,4,2})_{|_{\widetilde{\E}^{g-1}}}:\widetilde{\E}^{g-1}\rightarrow\E_2^{g-1}$.}
So a set of generators of  $\widetilde{\E}^{g-1}$ is given by the
generators of $\ker j$ and the preimages of the generators of
$\E_2^{g-1}$. By Remark \ref{ker}, $\ker j$ is generated by spins
of $Z_g$ and $Z'_g$ about appropriate loops not containing $P_1$
and $P_2$. So they correspond to slides of $Z_g$ and $Z'_g$ on
$\F_g$, which are elements of $\Y_g$ by Lemma \ref{slide}.
Moreover, by the induction hypothesis  $\E_2^{g-1}=\Y^{g-1}$.
Since we can suppose that the generators of $\Y^{g-1}$ keep $Z_g$
and $Z'_g$ fixed,  they correspond to  element of $\Y^{g}$. So, as
above, the isotopy class of $f$ is generated by
$\rho_{ij},\xi'_{ij},\xi_{ij},\theta'_{ij},\theta_{ij},\tau_i,\omega_i,\eta_{i1},\eta_{i2},\eta'_{i1},\eta'_{i2}$,
with $i,j=1,\ldots ,g-1$, $i\ne j$, and
$\xi'_{gj},\xi_{gj},\theta'_{gj},\theta_{gj},\eta_{g1},\eta_{g2},\eta'_{g1},\eta'_{g2},\tau_g$,
with $j=1,\ldots,g-1$. As a consequence
$\widehat{\E}_2^g\subseteq\mathcal{Y}^g$.
\end{proof}

\bigskip

\noindent{\bf Acknowledgements}

Work performed under the auspices of the G.N.S.A.G.A. of
I.N.d.A.M. (Italy) and the University of Bologna, funds for
selected research topics.

\bigskip\bigskip

\vspace{15 pt} {\noindent ALESSIA CATTABRIGA, Department of
Mathematics, University of Bologna, Piazza di Porta S. Donato, 5,
40126, Bologna (Italy). E-mail: cattabri@dm.unibo.it}

\vspace{15 pt} {\noindent MICHELE MULAZZANI, Department of
Mathematics, C.I.R.A.M., University of Bologna, Piazza di Porta S.
Donato, 5, 40126, Bologna (Italy). E-mail: mulazza@dm.unibo.it}

\end{document}